\documentclass[letter,leqno]{amsart}
\usepackage{amssymb}
\numberwithin{equation}{section}

\newtheoremstyle{theorem}{3pt}{3pt}%
{\it}
{}
{\bfseries}
{:}
{.5em}
{}
\theoremstyle{theorem}
\newtheorem{theorem}{Theorem}[section]

\newtheorem{proposition}[theorem]{Proposition}
\newtheorem{corollary}[theorem]{Corollary}

\newtheorem{definition}[theorem]{Definition}
\newtheorem{definitiontheorem}[theorem]{Theorem/Definition}

\newtheoremstyle{example}{3pt}{3pt}%
{}
{}
{\sc}
{:}
{.5em}
{}
\theoremstyle{example}
\newtheorem{example}[theorem]{Example}
\newtheoremstyle{remark}{3pt}{3pt}%
{}
{}
{\sc}
{:}
{.5em}
{}
\theoremstyle{remark}
\newtheorem{remark}{Remark}[section]
\newtheorem{remarks}{Remarks}[section]
\numberwithin{equation}{section}

\newcommand{\thismonth}{\ifcase\month\or
  January\or February\or March\or April\or May\or June\or
  July\or August\or September\or October\or November\or December\fi
  \space\number\year}
\makeatletter
\newcommand{\low}{\@ifnextchar^{}{^{\vphantom x}}}
\newcommand{\high}{\@ifnextchar_{}{_{\vphantom I}}}
\makeatother
\DeclareSymbolFont{script}{U}{eus}{m}{n}
\DeclareSymbolFontAlphabet{\mathscr}{script}
\DeclareMathSymbol{\EuWedge}{0}{script}{"5E}
\DeclareMathAlphabet{\mathrmsl}{OT1}{cmr}{m}{sl}
\newcommand{\rssymb}[2]{\newcommand{#1}{{\mathrmsl{#2}}}}
\newcommand{\calsymb}[2]{\newcommand{#1}{{\mathcal{#2}}}}
\newcommand{\bbsymb}[2]{\newcommand{#1}{{\mathbb{#2}}}}
\newcommand{\lieoper}[2]{\newcommand{#1}{\mathop
  {\mathfrak{#2}\null}\nolimits}}
\newcommand{\oper}[3][n]{\newcommand{#2}{\mathop
  {\mathrm{#3}\null}\ifx n#1\nolimits\else\limits\fi}}
\newcommand{\rsoper}[3][n]{\newcommand{#2}{\mathop
  {\mathrmsl{#3}\null}\ifx n#1\nolimits\else\limits\fi}}
\bbsymb\C{C} \bbsymb\F{F} \bbsymb\HQ{H}\bbsymb\N{N} \bbsymb\Q{Q}
\bbsymb\R{R} \bbsymb\U{U} \bbsymb\V{V} \bbsymb\W{W} \bbsymb\Z{Z}
\bbsymb\bbf{F} \bbsymb\bbk{K} \bbsymb\bbi{I} \bbsymb\bbl{L} \bbsymb\bbo{O}
\bbsymb\bbj{J}
\bbsymb\bby{Y}
\bbsymb\bbp{P}
\bbsymb\bba{A}
\calsymb\cA{A} \calsymb\cB{B} \calsymb\cC{C} \calsymb\cD{D} \calsymb\cE{E}
\calsymb\cF{F} \calsymb\cG{G} \calsymb\cH{H} \calsymb\cI{I} \calsymb\cJ{J}
\calsymb\cK{K} \calsymb\cL{L} \calsymb\cM{M} \calsymb\cN{N} \calsymb\cO{O}
\calsymb\cP{P} \calsymb\cQ{Q} \calsymb\cR{R} \calsymb\cS{S} \calsymb\cT{T}
\calsymb\cU{U} \calsymb\cV{V} \calsymb\cW{W} \calsymb\cX{X} \calsymb\cY{Y}
\calsymb\cZ{Z}
\renewcommand{\geq}{\geqslant} \renewcommand{\leq}{\leqslant}
\oper\End{End} \oper\Hom{Hom}                    
\oper\Sym{Sym} \oper\Skew{Skew}
\oper\Aut{Aut}                                   
\oper\GL{GL} \oper\SL{SL}\oper\Symp{Sp}
\oper\CO{CO} \oper\On{O} \oper\SO{SO} \oper\Pin{Pin} \oper\Spin{Spin}
\oper\CU{CU} \oper\Un{U} \oper\SU{SU}
\rsoper\Diff{Diff} \rsoper\SDiff{SDiff}
\lieoper\der{der}                                
\lieoper\gl{gl} \lieoper\sgl{sl}\lieoper\symp{sp}
\lieoper\co{co} \lieoper\so{so} \lieoper\spin{spin}
\lieoper\cu{cu} \lieoper\un{u}  \lieoper\su{su}
\rsoper\Vect{Vect} \rsoper\Ham{Ham}
\newcommand{\norm}[2][]{|\mkern-2mu|#2|\mkern-2mu|
  _{\lower1pt\hbox{${}_{#1}$}}}
\newcommand{\Norm}[2][]{\bigl|\mkern-3mu\bigr|#2\bigr|\mkern-3mu\bigr|
  _{\lower1pt\hbox{${}_{#1}$}}}

\newcommand{\Lie}{\cL}                      
\rsoper\dimn{dim}                           
\rsoper\grad{grad}                          
\rsoper\kernel{ker}\rsoper\image{im}        
\rsoper\alt{alt}   \rsoper\sym{sym}         
\rsoper\Ad{Ad}     \rsoper\ad{ad}           
\rsoper\CoAd{CoAd} \rsoper\coad{coad}       
\rsoper\trace{tr}  \rsoper\trfree{tf}       
\rsoper\detm{det}                           
\rsoper\Vol{Vol}                            
\rsoper\divg{div}                           
\rsoper\sign{sign}                          
\rssymb\iden{id}                            
\rssymb\vol{vol}                            
\oper\Imag{Im}\oper\Real{Re}                
\newcommand{\sd}{{\raise1pt\hbox{$\scriptscriptstyle +$}}}
\newcommand{\asd}{{\raise1pt\hbox{$\scriptscriptstyle -$}}}
\newcommand{\sdasd}{{\raise1pt\hbox{$\scriptscriptstyle\pm$}}}
\newcommand{\asdsd}{{\raise1pt\hbox{$\scriptscriptstyle\mp$}}}

\rsoper\scal{scal}
\def\kahl/{k\"ahler}
\def\Kahl/{K{\"a}hler}
\newcommand{\bfw}{{\mathbf w}}
\newcommand{\bfx}{{\mathbf x}}
\newcommand{\bfy}{{\mathbf y}}

\newcommand{\bfa}{{\mathbf a}}

\newcommand{\bfl}{{\mathbf l}}

\def\bba{{\mathbb A}}

\def\bbc{{\mathbb C}}

\def\bbf{{\mathbb F}}

\def\bbi{{\mathbb I}}
\def\bbj{{\mathbb J}}
\def\bbk{{\mathbb K}}
\def\bbl{{\mathbb L}}

\def\bbo{{\mathbb O}}
\def\bbp{{\mathbb P}}
\def\bbq{{\mathbb Q}}
\def\bbr{{\mathbb R}}

\def\bby{{\mathbb Y}}
\def\bbz{{\mathbb Z}}

\def\gra{\alpha}
\def\grb{\beta}

\def\grg{\gamma}
\def\gri{\iota}

\def\gro{\omega}

\def\grs{\sigma}
\def\grt{\tau}
\def\gru{\upsilon}

\def\grz{\zeta}

\def\grO{\Omega}

\def\cala{{\mathcal A}}
\def\calb{{\mathcal B}}
\def\calc{{\mathcal C}}
\def\calo{{\mathcal O}}

\def\cald{{\mathcal D}}

\def\calf{{\mathcal F}}

\def\calh{{\mathcal H}}

\def\calk{{\mathcal K}}

\def\calm{{\mathcal M}}
\def\caln{{\mathcal N}}
\def\calo{{\mathcal O}}

\def\calr{{\mathcal R}}
\def\cals{{\mathcal S}}

\def\calz{{\mathcal Z}}

\def\Se{Sasakian-Einstein }
\def\la#1{\hbox to #1pc{\leftarrowfill}}
\def\ra#1{\hbox to #1pc{\rightarrowfill}}

\def\gg{{\mathfrak g}}

\def\gl{{\mathfrak l}}

\def\gn{{\mathfrak n}}
\def\go{{\mathfrak o}}

\def\gt{{\mathfrak t}}
\def\gu{{\mathfrak u}}

\def\gA{{\mathfrak A}}

\def\gC{{\mathfrak C}}

\def\gF{{\mathfrak F}}

\def\gH{{\mathfrak H}}

\def\gT{{\mathfrak T}}

\def\gZ{{\mathfrak Z}}

\def\fract#1#2{\raise4pt\hbox{$ #1 \atop #2 $}}

\begin{document}

\title{Constructions in Sasakian Geometry}
\footnote{During the preparation of this work the first two authors
were partially supported by NSF grants DMS-0203219 and DMS-0504367.}

\author{Charles P. Boyer}
\address{Department of Mathematics and Statistics,
\newline \indent University of New Mexico,
\newline \indent Albuquerque, N.M. 87131}
\email{cboyer@math.unm.edu}

\author{Krzysztof Galicki}
\address{Department of Mathematics and Statistics,
\newline \indent University of New Mexico,
\newline \indent Albuquerque, N.M. 87131}
\email{galicki@math.unm.edu}

\author{Liviu Ornea}
\address{University of Bucharest, Faculty of Mathematics
\newline \indent
14 Academiei str., 70109 Bucharest, Romania}
\email{Liviu.Ornea@imar.ro}

\subjclass{} \keywords{Sasakian manifold, contact structures, join construction, contact
fiber bundles, toric contact manifolds}

\maketitle

\section{Introduction}

A Riemannian manifold $(M,g)$ is called a {\it Sasakian manifold}
if there exists a Killing vector field $\xi$ of unit length on $M$
so that the tensor field $\Phi$ of type $(1,1)$, defined by
$\Phi(X) ~=~ -\nabla_X \xi$, satisfies the condition $(\nabla_X
\Phi)(Y) ~=~ g(X,Y)\xi -g(\xi,Y)X$ fo any pair of vector fields
$X$ and $Y$ on $M.$ This is a curvature condition which can be
easily expressed in terms  the Riemann curvature tensor as
$R(X,\xi)Y ~=~ g(\xi,Y)X-g(X,Y)\xi$. Equivalently, the Riemannian
cone defined by $(\calc(M), \bar{g}, \Omega) ~=~ (\bbr_+\times M,
\ dr^2+r^2g, d(r^2\eta))$ is K\"ahler with the K\"ahler form
$\Omega=d(r^2\eta)$, where $\eta$ is the dual 1-form of $\xi$. The
4-tuple $\cals=(\xi,\eta,\Phi,g)$ is commonly called a {\it
Sasakian structure} on $M$ and $\xi$ is its {\it characteristic}
or {\it Reeb vector field}.

Sasakian geometry is a special kind of contact metric geometry
such that the structure transverse to the Reeb vector field $\xi$ is
K\"ahler and invariant under the flow of $\xi.$ In
fact $\eta$ is the contact 1-form, and $\Phi$ is a $(1,1)$ tensor
field which defines a complex structure on the contact subbundle
$\ker~\eta$ which annihilates $\xi.$
When both $M$ and the leaves of the foliation generated by $\xi$
are compact the Sasakian structure is called {\it quasi-regular},
and the space of leaves $X^{orb}$ is a compact K\"ahler orbifold.
In such a case $M$ is the total space of a circle {\it orbi-bundle}
(also called V-bundle) over $X^{orb}.$ Moreover, the 2-form
$d\eta$ pushes down to a K\"ahler form $\gro$ on $X^{orb}.$  Now
$\gro$ defines an integral class $[\gro]$ of the orbifold
cohomology group $H^2(X^{orb},\bbz)$ which generally is only a
rational class in the ordinary cohomology $H^2(X,\bbq).$

This construction can be inverted in the sense that given a
K\"ahler form $\gro$ on a compact complex orbifold $X^{orb}$ which
defines an element $[\gro]\in H^2(X^{orb},\bbz)$ one can construct
a circle  orbi-bundle on $X^{orb}$ whose orbifold first Chern
class is $[\gro].$ Then the total space $M$ of this orbi-bundle
has a natural Sasakian structure $(\xi,\eta,\Phi,g)$, where $\eta$
is a connection 1-form whose curvature is $\gro.$ The tensor field
$\Phi$ is obtained by lifting the almost complex structure $I$ on
$X^{orb}$ to the horizontal distribution $\ker~\eta$ and requiring
that $\Phi$ annihilates $\xi.$ Furthermore, the map $(M,g)\ra{1.3}
(X^{orb},h)$ is an orbifold Riemannian submersion. This is an
orbifold version of a well-known construction of Kobayashi. For
the essentials and more details on Sasakian geometry we refer the
reader to the forthcoming book \cite{BG05} of the first two authors.

The purpose of this paper is to describe in detail certain
constructions of new Sasakian manifolds from old ones. In Section
2 we generalize the join construction introduced by the first two
authors \cite{BG00a} in the case of quasi-regular \Se manifolds to
arbitrary quasi-regular Sasakian spaces. This construction is far
more flexible yielding a multitude of examples. Furthermore, owing
to the recent \Se metrics discovered on $S^5$ in \cite{BGK05}, we
are able to prove the existence of families of \Se metrics on
manifolds homeomorphic to $S^2\times S^5.$ However, determination
of the smooth structure is rather subtle and would ultimately
involve computation of the Kreck-Stolz invariants for these
manifolds \cite{KS88}.

In Section 3 we show how the join construction emerges as a
special case of Lerman's contact fibre bundle construction
\cite{Ler04} which under some additional assumptions can be
adapted to the Sasakian case. In particular, when both the base
and the fiber of the contact fiber bundle are toric we show that
the construction yields a new toric Sasakian manifold.

In the last section we study the toric Sasakian manifolds in
dimension 5. All compact, smooth, simply-connected, oriented
5-manifolds were classified by fundamental theorems of Smale and
Barden \cite{Sm62, Bar65}. In particular, the manifold (and its
unique smooth structure) is completely determined by
$H_2(M^5,\bbz)$ together with the second Stiefel-Whitney class map
$w_2:H_2(M,\bbz)\ra{1.2}\bbz_2$. The second Betti number $b_2(M)$,
the structure of the 2-torsion subgroup, and $w_2$ all provide
obstruction to the existence of various geometric structures on
such manifolds. For example, it is an elementary result that
torsion in the second homology group is the obstruction to the
existence of a free circle action on $M^5$. Moreover, vanishing of
the torsion is also a necessary and sufficient condition for the
existence of a regular contact structure as observed by Geiges
\cite{Gei91}. Remarkably, Oh proves the same condition is also a
necessary and sufficient for the existence of an effective $T^3$
action on $M^5$ \cite{Oh83} and Yamazaki shows that in such a case
one can always choose a $T^3$ action with a compatible toric
K-contact structure. We use these results to show that any
simply-connected compact oriented 5-manifold with vanishing
2-torsion admits a toric Sasakian structure. Furthermore, we prove
by explicitly constructing circle bundles over the blow-ups of
Hirzebruch surfaces that one can always find toric Sasakian
structures which are regular.

\medskip

{\bf Acknowledgements.} The authors thank Santiago Simanca and
Stefan Stolz for helpful discussions. L.O. is grateful to the
Department of Mathematics and Statistics of the University of New
Mexico in Albuquerque for warm hospitality and an excellent
research environment during the academic year 2004--2005. L.O.
also acknowledges partial funding from the Efroymson Foundation, as 
well as Grant 2-CEx-06-11-22/25.07.06.
K.G. would like to thank IHES and MPIM in Bonn for hospitality and
support. Finally, K.G. and L.O. would like to thank National
Research Council. In 2003/2004 our joint COBASE grant sponsored
short visits of L.O. at UNM and K.G. in Bucharest which visits
seeded the later collaboration.

\section{The Join Construction}\label{joinsection}

In this subsection we apply a construction due to Wang and Ziller \cite{WaZi90} to define
a
multiplication on the set of quasi-regular Sasakian orbifolds. This was done originally in
\cite{BG00a} in the case of Sasakian-Einstein orbifolds which is perhaps of more interest,
but there is an easy generalization to the strict Sasakian case. The idea is quite simple and
is
based on the fact that product of K\"ahler orbifolds is a K\"ahler orbifold in a natural way.

\begin{definition}\label{sewz.1}We denote by $\cals\calo$ the set of
compact quasi-regular Sasakian orbifolds, by $\cals\calm$ the
subset of $\cals\calo$ that are smooth manifolds, and by $\calr\subset
\cals\calm$ the subset of compact, simply connected, regular Sasakian manifolds.
The set $\cals\calo$ is topologized with the $C^{m,\gra}$ topology, and the
subsets are given the subspace topology.
\end{definition}

The set $\cals\calo$ is graded by dimension, that is,
$$\cals\calo =\bigoplus_{n=0}^\infty \cals\calo_{2n+1},$$
and similarly for $\cals\calm$ and $\calr.$ In the definition of Sasakian
structure it is implicitly assumed that $n>0.$ So we want to extend the
definition of a Sasakian structure to the case when $n=0.$ This can easily be done
since a connected one dimensional orbifold is just an interval with possible
boundary, or a circle.  So we can just take $\xi ={\partial\over \partial t},
\eta =dt,\Phi=0,$ with the flat metric $g=dt^2.$ In this case the space of
leaves $\calz$ of the characteristic foliation is just a point. The unit circle
$S^1$ with this structure will play the role of the identity in our monoid. Notice that the
identity is smooth.

For each pair of relatively prime positive integers $(k_1,k_2)$ we define a graded
multiplication
\begin{equation}\label{joinmaps}
\star_{k_1,k_2}:\cals\calo_{2n_1+1}\times \cals\calo_{2n_2+1}\ra{1.5}
\cals\calo_{2(n_1+n_2)+1}
\end{equation}
as follows: Let $\cals_1,\cals_2\in \cals\calo$ of dimension $2n_1+1$ and
$2n_2+1$ respectively. Since each orbifold $\cals_i$ has a quasi-regular Sasakian
structure,
its Reeb vector field generates a locally free circle action, and the quotient space by this
action has a natural orbifold structure $\calz_i$ \cite{Mol88}. Thus, there is a locally free
action of the 2-torus $T^2$ on the product orbifold $\cals_1\times \cals_2,$ and the
quotient
orbifold is the product of the orbifolds $\calz_i.$ (Locally free torus actions on orbifolds
have
been studied in \cite{HaSa91}). Now the Sasakian structure on $\cals_i$ determines a
K\"ahler structure $\gro_i$ on the orbifold $\calz_i$, but in order to obtain an integral
orbifold
cohomology class $[\gro_i]\in H^2(\calz_i,\bbz)$ we need to assure that the period of a
generic orbit is one. By a result of Wadsley \cite{Wad} the period function on a
quasi-regular
Sasakian orbifold is lower semi-continuous and constant on the dense open set of regular
orbits. This is because on a Sasakian orbifold all Reeb orbits are geodesics. Thus, by a
transverse homothety we can normalize the period function to be the constant $1$ on the
dense open set of regular orbits. In this case the K\"ahler forms $\gro_i$ define integer
orbifold cohomology classes $[\gro_i]\in H^2_{orb}(\calz_i,\bbz).$ If $Z_i$ denotes the
underlying
complex space associated with the orbifold $\calz_i$, one should not confuse
$H^*_{orb}(\calz_i,\bbz)$ with $H^*(Z_i,\bbz)$ nor with the Chen-Ruan cohomology of an
orbifold. $H^*_{orb}(\calz_i,\bbz)$ is the orbifold cohomology defined by Haefliger
\cite{Hae84}
(see also \cite{BG00a}). Notice, however, that we do have
$H^*_{orb}(\calz_i,\bbz)\otimes_\bbz\bbq \approx H^*(Z_i,\bbz)\otimes_\bbz\bbq.$ Now
each
pair of positive integers $k_1,k_2$ give a K\"ahler form $k_1\gro_1+k_2\gro_2$ on the
product. Furthermore,  $[k_1\gro_1+k_2\gro_2]\in H^2_{orb}(\calz_1\times \calz_2,\bbz),$
and
thus
defines an $S^1$ V-bundle over the orbifold $\calz_1\times \calz_2$ whose total space is
an
orbifold that we denote by $\cals_1\star_{k_1,k_2} \cals_2$ and refer to as the
$(k_1,k_2)$-{\it join} of $\cals_1$ and $\cals_2.$ Furthermore, $\cals_1\star_{k_1,k_2}
\cals_2$ admits a quasi-regular Sasakian structure \cite{BG05} by choosing a connection
1-form on $\cals_1\star_{k_1,k_2} \cals_2$ whose curvature is
$\pi^*(k_1\gro_1+k_2\gro_2).$ This Sasakian structure is unique up to a gauge
transformation of the form $\eta\mapsto \eta +d\psi$ where $\psi$ is a smooth basic
function. This defines the maps in \ref{joinmaps}. If $\cals_i$ are quasi-regular Sasakian
structures on the compact manifolds $M_i,$ respectively, we shall use the notation
$\cals_1\star_{k_1,k_2} \cals_2$ and $M_1\star_{k_1,k_2} M_2$ interchangeably
depending
on whether we want to emphasize the Sasakian or manifold nature of the join. Notice also
that if $\gcd(k_1,k_2)=m$ and we define $(k'_1,k'_2)=(\frac{k_1}{m},\frac{k_2}{m}),$ then
$\gcd(k'_1,k'_2)=1$ and $M_1\star_{k_1,k_2} M_2\approx (M_1\star_{k'_1,k'_2}
M_2)/\bbz_m.$ In this case the cohomology class $k'_1\gro_1+k'_2\gro_2$ is indivisible in
$H^2_{orb}(\calz_1\times \calz_2,\bbz).$ Note also that $M_1\star_{k_1,k_2} M_2$ can be
realized as the quotient space $(M_1\times M_2)/S^1(k_1,k_2)$ where the $S^1$   action
is
given by the map
\begin{equation}\label{s1action}
(x,y)\mapsto (e^{ik_2\theta}x,e^{-ik_1\theta}y).
\end{equation}

We are interested in restricting the map $\star_{k_1,k_2}$ of \ref{joinmaps} to the subset
of
smooth Sasakian manifolds, that is in the map
\begin{equation}\label{joinmaps2}
\star_{k_1,k_2}:\cals\calm_{2n_1+1}\times \cals\calm_{2n_2+1}\ra{1.5}
\cals\calo_{2(n_1+n_2)+1}.
\end{equation}
If $M_1$ and $M_2$ are quasi-regular Sasakian manifolds,
we are interested under what conditions the orbifold $M_1\star_{k_1,k_2}M_2$ is a
smooth
manifold. Let $\upsilon_i$ denote the order of the quasi-regular Sasakian manifolds
$M_i,$ that is, $\upsilon_i$ is the lcm of the orders of the leaf holonomy groups of $M_i.$
The
following proposition is essentially Proposition 4.1 of \cite{BG00a}:

\begin{proposition}\label{smooth}
For each pair of relatively prime positive integers $k_1,k_2$, the orbifold
$M_1\star_{k_1,k_2}M_2$ is
a smooth quasi-regular Sasakian manifold if and only if
$\gcd(\upsilon_1k_2,\upsilon_2k_1)=1.$ In particular, if $M_i$ are regular Sasakian
manifolds,
then so is $M_1\star_{k_1,k_2}M_2.$
\end{proposition}

Generally, given two known Sasakian manifolds $M_1$ and $M_2,$ it can be quite difficult
to
compute the diffeomorphism type of $M_1\star_{k_1,k_2}M_2,$ but some information can
be
obtained. For example, we have

\begin{proposition}\label{assocVbunjoin}
Let $M_1$ and $M_2$ be compact quasi-regular Sasakian manifolds
and assume that $\gcd(\upsilon_1k_2,\upsilon_2k_1)=1.$  Then
$M_1\star_{k_1,k_2}M_2$ is the associated $S^1$ orbibundle over
$\calz_1$ with fibre $M_2/\bbz_{k_2}.$ In particular, if $M_1$ is
regular and $k_2=1$ then for each positive integer $k,$
$M_1\star_{k,1}M_2$ is an $M_2$-bundle over the K\"ahler manifold
$Z_1.$
\end{proposition}

\begin{proof}
Following \cite{WaZi90} we break up the $S^1$ action on $M_1\times M_2$ into stages.
First
divide by the subgroup
$\bbz_{k_2}$ of the circle group $S^1(k_1,k_2)$ defined by Equation \ref{s1action} giving
$M_1 \times (M_2/\bbz_{k_2}).$ Letting $[y]\in M_2/\bbz_{k_2}$ denote the equivalence
class
of $y\in M_2,$ we see that the quotient group
$S^1/\bbz_{k_2}$ acts on $M_1 \times M_2/\bbz_{k_2}$ by $(x,[y])\mapsto
(e^{i\theta}x,[e^{-i\frac{k_1}{k_2}\theta}y])$ which identifies $M_1\star_{k_1,k_2}M_2$ as
the
orbibundle over $\calz_1$ with fibre $M_2/\bbz_{k_2}$ associated to the principal $S^1$
orbibundle $\pi_1:M_1\ra{1.5} \calz_1.$ This
proves the result.
\end{proof}

Recall \cite{BG05} the {\it type} of a Sasakian structure. A Sasakian structure
$(\xi,\eta,\Phi,g)$ is said to be of {\it positive (negative) type} if the first Chern class
$c_1(\calf_\xi)$ of the characteristic foliation is represented by a positive (negative)
definite
$(1,1)$-form. If either of these two conditions
is satisfied $(\xi,\eta,\Phi,g)$ is said to be of {\it definite type}, and otherwise
$(\xi,\eta,\Phi,g)$ is of {\it indefinite type}. $(\xi,\eta,\Phi,g)$ is said to be of {\it null type} if
$c_1(\calf_\xi)=0.$ We often just say `a positive Sasakian structure' instead of `a Sasakian
structure of positive type', etc. It will also be convenient to write $c_1(\cals)$ instead of
$c_1(\calf_\xi)$ even though $c_1$ is independent of the Sasakian structure in the
deformation class $\gF(\xi)$ \cite{BG05}.

\begin{proposition}\label{typeprop}
The $(k_1,k_2)$-join of two positive, negative, or null compact quasi-regular Sasakian
manifolds is positive, negative, or null, respectively.
\end{proposition}

\begin{proof}
This follows from the fact that for any quasi-regular Sasakian structure on a compact
manifold $M$ we have an orbifold submersion $\pi:M\ra{1.5} \calz$ satisfying
$c_1(\calf_\xi)=\pi^*c_1^{orb}(\calz)$ as real cohomology classes. So the sign or
vanishing of
$c_1(\calf_\xi)$ and $c_1^{orb}(\calz)$ coincide. Furthermore, for any pair of integers
$(k_1,k_2)$ we have
$$c_1(\cals_1\star_{k_1,k_2} \cals_2)=\pi^*c_1^{orb}(\calz_1\times \calz_2)=
\pi^*c_1^{orb}(\calz_1) +\pi^*c_1^{orb}(\calz_2)=c_1(\cals_1)+c_1(\cals_2).$$
Now suppose that the Sasakian structures $\cals_1$ and $\cals_2$ are both definite of the
same type. Then $c_1(\cals_1)+c_1(\cals_2)$ can be represented by either a positive
definite
or negative definite basic $(1,1)$ form. The null case is clear.
\end{proof}

Next we give examples of Wang and Ziller \cite{WaZi90} where the topology can be
ascertained,
even the homeomorphism and diffeomorphism type in certain cases.

\begin{example}\label{nontrivialsks2} The Wang-Ziller manifolds:\index{Wang-Ziller
manifolds}
Let $M^{p_1,p_2}_{k_1,k_2}$ denote $S^{2p_1+1}\star_{k_1,k_2}S^{2p_2+1}.$ This is the
$S^1$-bundle over $\bbc\bbp^{p_1}\times \bbc\bbp^{p_2}$ whose first
Chern class is $k_1[\gro_1]+k_2[\gro_2]$ where $\gro_i$ is the standard K\"ahler class of
$H^*(\bbc\bbp^{p_i},\bbz)$  and $k_i\in \bbz^+.$ By Proposition
\ref{smooth} $M^{p_1,p_2}_{k_1,k_2}$ admits regular Sasakian structures, and by
Proposition
\ref{typeprop} they are
positive. Furthermore, if $\gcd(k_1,k_2)=1$ (which we assume hereafter) the manifolds
$M^{p_1,p_2}_{k_1,k_2}$
are simply connected.
To analyze the manifolds $M^{p_1,p_2}_{k_1,k_2}$ we follow Wang and Ziller
\cite{WaZi90}
and
consider the free $T^2$ action on $S^{2p_1+1}\times S^{2p_2+1}$ defined by
$(\bfx,\bfy)\mapsto
(e^{i\theta_1}\bfx,e^{i\theta_2}\bfy)$ where $(\bfx,\bfy)\in \bbc^{p_1+1}\times
\bbc^{p_2+1}.$
The
quotient space is $\bbc\bbp^{p_1}\times \bbc\bbp^{p_2},$ and $M^{p_1,p_2}_{k_1,k_2}$
can
be identified with the
quotient of $S^{2p_1+1}\times S^{2p_2+1}$ by the circle defined by
$(\bfx,\bfy)\mapsto (e^{ik_2\theta}\bfx,e^{-ik_1\theta}\bfy).$ Now the free part of
$H^2(M^{p_1,p_2}_{k_1,k_2},\bbz)$ has a single generator $\grg,$ and  letting $\pi:
M^{p_1,p_2}_{k_1,k_2}\ra{1.3}
\bbc\bbp^{p_1}\times \bbc\bbp^{p_2}$ denote the natural bundle projection, we see that
the
classes
$[\gro_i]$ pull back as $\gri_*\pi^*[\gro_1]=k_2\grg,\gri_*\pi^*[\gro_2]=-k_1\grg.$ Here by
abuse of notation
we
let $\pi^*[\gro_i]$ also denote the basic classes in $H^2_B(\calf_\xi).$
Furthermore, the basic first Chern class is
$c_1(\calf_\xi)=(p_1+1)\pi^*[\gro_1]+(p_2+1)\pi^*[\gro_2]$, so we get
$$c_1(\cald)=\gri_*c_1(\calf_\xi)=(p_1+1)\gri_*\pi^*[\gro_1]+(p_2+1)\gri_*\pi^*[\gro_2]
=(k_2(p_1+1)-k_1(p_2+1))\grg.$$
Thus, we have
\begin{equation}\label{spinformula}
w_2(M^{p_1,p_2}_{k_1,k_2})= (k_2(p_1+1)+k_1(p_2+1))\grg \mod 2.
\end{equation}
In certain cases one can determine the manifold completely
\cite{WaZi90}. For example, consider $p_1=k_2=1,p_2=q,k_1=k$ in
which case $M^{1,q}_{k,1}$ is an $S^{2q+1}$-bundle over $S^2.$ The
$S^k$-bundles over $S^2$ are classified by $\pi_1(SO(k+1))\approx
\bbz_2$ \cite{Stee51}. So there are precisely two
$S^{2q+1}$-bundles over $S^2,$ and they are distinguished by
$w_2.$ From Equation \eqref{spinformula} we get
$w_2(M^{1,q}_{k,1})= k(q+1)\grg \mod 2.$ Thus, if $q$ is odd or
$k$ is even, we get the trivial bundle $S^2\times S^{2q+1};$
whereas, if $q$ is even and $k$ is odd, we get the unique
non-trivial $S^{2q+1}$-bundle over $S^2.$ This gives an infinite
number of distinct deformation classes of regular positive
Sasakian structures on these manifolds. In dimension five
$(p_1=p_2=1)$ we can do somewhat better. In fact for any pair of
relatively prime positive integers $(k_1,k_2)$ $M^{1,1}_{k_1,k_2}$
is diffeomorphic to $S^2\times S^3;$ whereas later in Section we
construct a positive Sasakian structure on the non-trivial
$S^3$-bundle over $S^2$ as well as a family of indefinite Sasakian
structures. Notice that in this case $c_1(\cald)=2(k_1-k_2)\grg.$
\end{example}

Summarizing from this example gives

\begin{corollary}\label{WZcor}
The manifolds $M^{p,q}_{k_1,k_2}$ all admit Sasakian metrics with positive Ricci curvature.
In particular, the manifolds $S^2\times S^{2q+1}$ as well as the non-trivial
$S^{2q+1}$-bundle over $S^2$ admit Sasakian metrics of positive Ricci curvature.
\end{corollary}

Wang and Ziller \cite{WaZi90} were able to prove the existence of positive Einstein
metrics on
these manifolds. As we have shown these manifolds always admit positive Sasakian
metrics,
but the Einstein metrics on $S^2\times S^{2q+1}$ are not generally Sasakian-Einstein. For
example, to get a \Se metric a particular join is necessary \cite{BG00a}. The Wang-Ziller
construction gives \Se metrics on $S^3\star_{2,q+1}S^{2q+1}$ for $q$ even, and
$S^3\star_{1,\frac{q+1}{2}}S^{2q+1}$ for $q$ odd. These are all non-trivial fibre bundles
over
$S^2$ whose fibres are the appropriate lens spaces.

We are interested in when the join of two Sasakian $\eta$-Einstein
manifolds of the same type is a Sasakian $\eta$-Einstein manifold
of that type. In the case of \Se structures this was described in
\cite{BG00a}.  In that case we had to choose the pair $(k_1,k_2)$
to be the relative Fano indices of the two Sasakian manifolds. The
concept of index applies equally well to the negative definite
case, but the name Fano is inappropriate.  So assuming that
$c_1^{orb}(\calz)$ is either positive or negative definite, we
write $|c_1^{orb}(\calz)|$ to mean $c_1^{orb}(\calz)$ in the
positive case and $-c_1^{orb}(\calz)$ in the negative case. Then
we can  define the {\it divisibility index} or just {\it index} of
the polarized K\"ahler orbifold $(\calz,|c_1^{orb}(\calz)|)$ to be
the largest integer $I=I(\calz)$ such that $c_1^{orb}(\calz)/I$
defines an integral class in $H^2_{orb}(\calz,\bbz).$ We also say
that the polarized K\"ahler orbifold $(\calz,|c_1^{orb}(\calz)|)$
is {\it indivisible}. Thus, if $c_1^{orb}(\calz)$ is represented
by either a positive definite or negative definite $(1,1)$-form,
there is a unique indivisible polarized K\"ahler orbifold
$(\calz,|c_1^{orb}(\calz)|/I)$ that gives rise to a quasi-regular
Sasakian orbifold structure $\cals.$ So we can consider the index
$I$ to be an invariant of the Sasakian structure $\cals$ and write
$I=I(\cals)$ as well. In this regard, we also say that a
quasi-regular definite Sasakian structure $\cals$ on $M$ is {\it
indivisible}. By Theorem 2.1 of \cite{BG00a} indivisible positive
quasi-regular Sasakian manifolds are simply connected. For $i=1,2$
we define the {\it relative indices} of a pair of quasi-regular
definite Sasakian structures $(\cals_1,\cals_2)$ of the same type
by
$$l_i= \frac{I(\cals_i)}{\gcd(I(\cals_1),I(\cals_2))}.$$
Then, $\gcd(l_1,l_2)=1,$ and if $\gcd(\upsilon_1l_2,\upsilon_2l_1)=1,$ then
$M_1\star_{l_1,l_2}M_2$ is a smooth indivisible Sasakian manifold. If in addition $M_1$
and
$M_2$ are both positive Sasakian manifolds, $M_1\star_{l_1,l_2}M_2$ is a simply
connected
positive Sasakian manifold.

We now consider the join of two Sasakian $\eta$-Einstein manifolds. The case of positive
$\eta$-Einstein manifolds is essentially the same as for \Se manifolds. Of course, in the
positive case there are obstructions to satisfying the Monge-Amp\`ere equation on a
compact
K\"ahler orbifold or the transverse
Monge-Amp\`ere equation on a Sasakian manifold; whereas, in both
the negative case and null cases, there are no such obstructions. So for any null Sasakian
structure or negative Sasakian structure such that $c_1(\calf_\xi)$ is a multiple of
$[d\eta]_B$
on a compact manifold there exists a compatible Sasakian $\eta$-Einstein metric.

\begin{proposition}\label{etaEjoin}
Let $l_i$ be the relative indices for a pair of quasi-regular definite Sasakian
$\eta$-Einstein manifolds $M_i$ of the same type, respectively. Then the join
$M_1\star_{l_1,l_2}M_2$ admits a quasi-regular Sasakian $\eta$-Einstein structure of that
type.
\end{proposition}

We can now obtain new \Se metrics by combining Proposition
\ref{assocVbunjoin} with the results in \cite{BGK05,BGKT05,GhKo05}. For example take
$M_1=S^3$with its canonical round sphere Sasakian structure, and $M_2=S^5_\bfw$ with
one of the 68 deformation classes of \Se structures on $S^5$ found in \cite{BGK05} or
one of
the 12 \Se structures in \cite{GhKo05}. Let $L^5(l_2)$ denote the lens space
$S^5/C_{l_2}$
where $C_{l_2}\approx \bbz_{l_2}$ is the cyclic subgroup of the circle group $S^1_\bfw$
generated by the Reeb vector field of the corresponding Sasakian structure, and $l_1$ is
the
reduced Fano index of $S^5_\bfw$ with respect to $S^3.$ It is straightforward to
compute
$l_2=l_2(\bfw)$ as a function of $\bfw.$ Then we have

\begin{theorem}\label{s2s5}
If $\gcd(l_2,\gru_2)=1$ then $S^3\star_{l_1,l_2}S^5_\bfw$ is  the
total space of the fibre bundle over $S^2$ with fibre the lens
space $L^5(l_2),$ and it admits \Se metrics.
In particular, for 16 different weight vectors $\bfw$, the manifold
$S^3\star_{2,1}S^5_\bfw$ is
homeomorphic to $S^2\times S^5$ and admits Sasakian-Einstein
metrics including one 10-dimensional family. Moreover, for 3 different weight vectors
$\bfw$,
the manifold $S^3\star_{1,1}S^5_\bfw$ is homeomorphic to $S^2\times S^5$ and admits
Sasakian-Einstein metrics.
\end{theorem}

\begin{proof}
As mentioned above this is constructed using Proposition \ref{assocVbunjoin} with
$M_1=S^3$ with its standard
round sphere \Se structure, and $M_2$ one of the \Se structures on $S^5$ mentioned
above.
So the first statement follows. To prove the second statement we need to compute the
relative indices. Since $I(S^3)=2$ for the standard Sasakian structure on $S^3$, we need
to
consider two cases, namely when $I(S^5_\bfw)=1,$ or $2.$ In both cases the relative
index
$l_2=1,$ so we have $S^3\star_{l_1,1}S^5_\bfw$ which is an $S^5$ bundle over $S^2.$
These are classified by their second Stiefel-Whitney class $w_2.$ But by construction
the
orbifold first Chern class of $\bbc\bbp^1\times \calz$ is proportional to the first Chern
class of
the $S^1$ orbibundle defining $S^3\star_{l_1,1}S^5_\bfw.$ This implies that
$w_2(S^3\star_{l_1,1}S^5_\bfw)$ vanishes \cite{BGN03a}.  Now all of the Sasakian
structures on $S^5$ can be represented as links of Brieskorn-Pham polynomials of the
form
$f=z_0^{a_0}+z_1^{a_1}+z_2^{a_2}+z_3^{a_3}$ or deformations thereof. So we need to
compute the Fano index $I$ for the 68 Brieskorn polynomials representing $S^5$ that
admit
\Se metrics found in \cite{BGK05}, and the 12 found in \cite{GhKo05}. Now it is easy to
see
that in terms of the Brieskorn exponents the Fano index takes the form
$$I={\rm lcm}(a_0,a_1,a_2,a_3)\bigl(\sum_{i=0}^3\frac{1}{a_i}-1\bigr).$$
It is easy to write
a Maple program to determine the Fano index of the 80 cases.
There are 16 with $I=1$ and 3 with $I=2$ giving 19 in all. For example the 10 parameter
family given by Example 41 in \cite{BGK05} has $\bfa=(2,3,7,35),$ and one easily sees that
$I=1.$
\end{proof}

\begin{remarks} We remark that only the $S^5_\bfw$ found in \cite{BGK05} give rise to
\Se
metrics on $S^2\times S^5;$ the ones found in \cite{GhKo05} all have Fano index greater
than $2$ and give non-trivial lens spaces. In fact the largest $l_2$ obtained is $89.$
Generally, assuming $\gcd(l_2,\gru_2)=1$, an easy spectral sequence argument shows
that
the manifolds $S^3\star_{l_1,l_2}S^5_\bfw$ are simply connected with the rational
homology
type of $S^2\times S^5,$ but with $H^4(S^3\star_{l_1,l_2}S^5_\bfw,\bbz)\approx
\bbz_{l_2}.$
\end{remarks}

It is now quite straightforward to apply Proposition
\ref{assocVbunjoin} to many other cases. For example, we can
consider the join $S^3\star_{2,1} k(S^2\times S^3)_\bfw$ or
$S^3\star_{1,1}k(S^2\times S^3)_\bfw$ where $k(S^2\times
S^3)_\bfw$ is any of the \Se manifolds consider in
\cite{BGN03c,BGN02b,BG03} with Fano index $I=1$ in the first case,
and $I=2$ in the second. This gives \Se metrics on manifolds
whose rational cohomology can be determined as in \cite{BG00a}.
The higher index cases in \cite{BGN03c} can also be treated as long as the relative index $l_2$ is 
relatively
prime to the order of the Sasakian structure of $k(S^2\times
S^3)_\bfw.$

Recall (cf. \cite{BG05,BGM06}) that the real Heisenberg group
$\gH_{2n+1}(\bbr)$ admits a homogeneous Sasakian structure with its standard 1-form
$\eta= dz-\sum_iy^idx^i.$ As a manifold
$\gH_{2n+1}(\bbr)$ is just $\bbr^{2n+1}$ which can be
realized in terms of $n+2$ by $n+2$ nilpotent matrices of the form
\begin{equation}\label{Heisenbergmatrix}
 \left(
\begin{matrix}1 &x_1 &\cdots  &x_n& z \\
                    0 &1 &0 &\cdots & y_1\\
                \vdots && \ddots & \cdots &  \vdots \\
                   0 & \cdots & 0&1 & y_n \\
                   0 &\cdots &0 & 0& 1
       \end{matrix}
    \right).
\end{equation}
If we consider the discrete subgroup $\gH_{2n+1}(\bbz)$ of $\gH_{2n+1}(\bbr)$ defined
by
the matrices \ref{Heisenbergmatrix} with integer entrees, the quotient manifold
$\caln_{2n+1}
=\gH_{2n+1}(\bbr)/\gH_{2n+1}(\bbz)$ is a nilmanifold with an induced Sasakian structure.
As a coset space it is also a homogeneous manifold, but the homogeneous structure and
Sasakian structure are incompatible. The Reeb vector field generates the one dimensional
center $\gZ(\gH_{2n+1}(\bbr))$ of the group $\gH_{2n+1}(\bbr),$ and thus
$\gZ(\gH_{2n+1}(\bbr))$ induces an $S^1$ action on the quotient space $\caln_{2n+1}.$
This
$S^1$ is the connected component of the group $\gA\gu\gt$ of Sasakian automorphisms
of
$\caln_{2n+1}$ and makes $\caln_{2n+1}$ the total space of an $S^1$ bundle over the
principally polarized Abelian variety $\cala =T^{2n}$ with its standard complex structure.
One
can obtain many regular Sasakian structures on
$\caln_{2n+1}$ by considering it as a circle bundle over a polarized Abelian variety
and deforming the complex structure. Now we can form the join
$\caln_{2n+1}\star_{k_1,k_2} M$ where $M$ is a regular Sasakian manifold which by
Proposition \ref{assocVbunjoin} can be thought of as an $M/\bbz_{k_2}$-bundle over
$T^{2n}.$

Recall that a locally conformal K\"ahler manifold is a complex manifold which admits a
covering endowed with a K\"ahler metric with respect to which the group of deck transformations
acts by holomorphic homotheties
(cf. \cite{DrOr98}). The subclass of Vaisman manifolds can be characterized in terms of
Sasakian
geometry as follows (cf. \cite{OrVe03}): Any compact Vaisman manifold $P$ is a
suspension over a
circle, with fibre a Sasakian manifold $M$. Moreover, there exist a Sasakian automorphism
$\varphi$ of
$M$ and a positive $q$ such that $P$ is isomorphic with the quotient of the Riemannian
cone
$(M\times
\mathbb{R}_+, t^2g_M+dt^2)$ by the cyclic group generated  by $(x, t)\mapsto
(\varphi(x),qt)$.It is clear from the definition that the product of
two l.c.K. structures is not, in general, l.c.K.
Moreover, the product of two Vaisman manifolds might not be Vaisman: \emph{e.g.} the
product of two
Hopf surfaces has $b_1=2$ which prevents it to admit a Vaisman structure, for which
$b_1$ should be odd. Instead, we can combine the join construction with the structure theorem to 
define a join
of quasi-regular compact Vaisman manifolds. Note that, in fact, this is not restrictive, since any
compact Vaisman structure can be deformed to a quasi-regular one \cite{OrVe05}.
Summing up, we have:
\begin{definition}
Let $P_1, P_2$ be two compact, quasi-regular Vaisman manifolds and let $M_1, M_2$ be
the respective
Sasakian manifolds provided by the structure theorem. Then, for each $k_1, k_2\in
\mathbb{Z}$, the
suspension over the circle with fibre $M_1\star_{k_1,k_2}M_2$ is the {\bf join} of $P_1$ and
$P_2$.
\end{definition}

\medskip
\section{Contact Fibre Bundles and Toric Sasakian Structures}\label{contfib}
\medskip

The aim of this section is to briefly discuss a construction due to Lerman \cite{Ler04b} that
allows one to construct K-contact structures on the total space of a fibre bundle whose
fibres are K-contact. This construction generalizes the join construction described in Section 
\ref{joinsection} as well
as the fibre join construction of Yamazaki \cite{Yam99}.
Actually one can work within the pure contact setting, and it is the contact
analog of symplectic fibre bundles described in \cite{GuLeSt96}. Recall that a {\it contact structure}
on an oriented
manifold $M$ is an equivalence class of 1-forms $\eta$ which satisfy $\eta\wedge (d\eta)^n\neq 
0,$
where two such
forms $\eta,\eta'$ are equivalent if there is a nowhere vanishing smooth function $f$ such that
$\eta'=f\eta.$
Alternatively, a contact structure is a maximally non-integrable codimension one subbundle $\cald$
of the tangent
bundle $TM.$ The relation between the two descriptions is $\cald=\ker~\eta.$ We denote a contact
manifold by the
pair $(M,\cald)$, and let $\gC\go\gn(M,\cald)$ denote the group of contactomorphisms of
$(M,\cald),$ that is, the
subgroup of the group of diffeomorphisms of $M$ leaving the contact bundle $\cald$ invariant. The
contact manifold
$(M,\cald)$ is said to be {\it co-oriented} if the subbundle $\cald$ is oriented. The subgroup of
$\gC\go\gn(M,\cald)$
which fixes the orientation is denoted by $\gC\go\gn(M,\cald)^+.$
Here is Lerman's definition of a contact fibre bundle.

\begin{definition}\label{contfibdef}
A fibre bundle $F\ra{1.5}M\fract{\pi}{\ra{1.5}} B$ is called a {\bf contact fibre bundle} if
\begin{enumerate}
\item $F$ is a co-oriented contact manifold with contact bundle $\cald.$
\item There are an open cover $\{U_i\}$ of $B$ and local trivializations
$\phi_i:\pi^{-1}(U_i)\ra{1.5} U_i\times F$ such that for every point $p\in U_i\cap U_j$ the
transition functions $\phi_j\circ \phi^{-1}_i|_{\{p\}\times F}$ are elements of
$\gC\go\gn(F,\cald)^+.$
\end{enumerate}
\end{definition}

We need the notion of fatness of a bundle due to Weinstein \cite{Wei80}. Let $\gra$ be  a
connection 1-form  in a principal bundle $P(M,G)$ with Lie group $G,$ and let
$\grO=D\gra$
denote its curvature 2-form. Let $S\subset \gg^*$ be any subset in the dual $\gg^*$ of
the Lie algebra $\gg$ of $G.$ We say that the connection $\gra$  is {\it fat on $S$} if the
bilinear map
\begin{equation}\label{fateqn}
\mu\circ \grO:\calh P\times \calh P: \ra{2.5} \bbr
\end{equation}
is non-degenerate for all $\mu\in S.$ In
particular, if $G$ is a torus, the bundle $\pi:P\rightarrow B$ is
identified, up to a gauge transformation, by a connection form $A$ such that
$dA=\pi^*\omega$ with $[\omega]\in H^2(B,\bbz)$. Then, if $\omega$ is non-degenerate, that is a
symplectic form, $A$ is certainly fat on the image of the moment map.

We first recall the main lines of the construction, not in full generality,
but adapted to our needs. Let $\pi:P\rightarrow B$ be a principal $G$-bundle endowed with a
connection
$A$ (we don't distinguish between the connection and its 1-form). Let $F$ be
a K-contact manifold, with fixed contact form $\eta_F$ and Reeb field $\xi_F$.
Suppose $G\subset \mathrm{Aut}(F,\eta)$, \emph{i.e.} it acts (from the left)
on $F$ by strong contactomorphisms and denote by $\Psi:F\rightarrow
\mathfrak{g}^*$ the associated momentum map. Then Lerman proves:

\begin{theorem}\cite{Ler04b}
In the above setting, if the connection $A$ is fat at all the points of the
image of the momentum map $\Psi$, then the total space $M$ of the associated
bundle $P\times_GF$ admits a $K$-contact structure.
\end{theorem}

We are interested in the case that the underlying almost CR
structure of the K-contact structure is integrable. In this case
the manifold $P\times_GF$ will be Sasakian. One way of guarantying
this occurs in the toric setting, so we now give a brief review of
toric contact geometry. This was begun in \cite{BM93}, continued
in \cite{BG00b}, and completed in \cite{Ler02a}. Let $(M,\cald)$
be a co-oriented contact manifold,

\begin{definition}\label{toriccontact}
A {\bf toric contact manifold} is a triple $(M,\cald,T^{n+1})$ where $(M,\cald)$ is an
a co-oriented contact manifold of dimension $2n+1$ with an effective action
$$\cala:T^{n+1}\ra{1.5} \gC\go\gn(M,\cald)^+$$
of a $(n+1)$-torus $T^{n+1}.$
\end{definition}

In \cite{BG00b} the first two authors introduced the notion of a contact toric structure of Reeb type.
\begin{definition}\label{Reebtype}
We say that a torus action $\cala:T^{n+1}\ra{1.3} \gC\go\gn(M,\cald)^+$ is of {\bf Reeb type}
if there are a contact 1-form $\eta$ of the contact structure $\cald$ and an element
$\varsigma\in \gg$ such that $X^\varsigma$ is the Reeb vector field of $\eta.$
\end{definition}

Fixing a contact form $\eta$ it is easy to see that the action $\cala$ of
$T^{n+1}$ is of Reeb type if and only if there is an
element $\grt$ in the Lie algebra $\gt_{n+1}$ of $T^{n+1}$ such that
$\eta(X^\grt) >0.$ Note that when $T^{n+1}$ acts properly we can always fix a contact 1-form 
$\eta$
without loss of
generality by using a slice theorem. In this case the relevant group is the subgroup
$\gC\go\gn(M,\eta)\subset
\gC\go\gn(M,\cald)^+$ of contactomorphisms leaving $\eta$ invariant. We are now ready for

\begin{theorem}\label{Sasfibre}
Let $F^{2n+1}$ be a compact toric contact manifold of Reeb type and with
torus
$T^{n+1}\subset \gC\go\gn(F,\eta)$. Let $\pi:P\rightarrow B$ be a
principal $T^{n+1}$
bundle over a toric compact symplectic manifold $B$.
Then $P\times_{T^{n+1}}F$ is a toric Sasakian manifold.
\end{theorem}

\begin{proof}
Choosing a connection $A$ as above whose curvature is non-degenerate, this will be fat and Lerman's
construction applies.
From \cite{BG00b}, $F$ has a compatible Sasakian structure, in particular
it is $K$-contact.
Note that this Sasakian structure is toric. Then by  \cite{Ler04b}, in our
hypothesis, $P\times_{T^{n+1}}F$ has a $K$-contact structure. We claim that
this is toric, of Reeb type, and hence the result follows by applying again \cite{BG00b}.
To prove our claim, we show that:
\begin{enumerate}
\item The Hamiltonian action of  $T^m$  on $B$ ($\dim B=2m$)  lifts to a $T^m$
action on $P$.
\item This lifted action extends to  $P\times_{T^{n+1}}F$ preserving the
contact form.
\item The action of $T^m$ on $P\times_{T^{n+1}}F$ commutes with the action of
$T^{n+1}$,
hence, as it leaves the contact form $\eta$ invariant, it induces an action
of $T^{n+m+1}$ on
$P\times_{T^{n+1}}F$. Then we only need to see that the Reeb field of
$P\times_{T^{n+1}}F$
is generated by the $T^{n+1}$ action.
\end{enumerate}
To prove (i), it is enough to show that $T^m$ lifts to an action that
preserves the fat connection
$A$. Denote $\{Y_i\}$ the generators of the $T^m$ action on $B$. We need to
construct lifts
$\tilde Y_i$ on $P$ such that  $\Lie_{\tilde Y_i}A=0$. Define them as
$$\tilde Y_i=\hat Y_i+a^\alpha_iX_\alpha,$$
where the $\hat{}$ refers to horizontal lifts, $X_\alpha$ are vertical
fields and the functions
$a_i^\alpha$ need to be determined. If we let $\{e_\alpha\}$ be a basis of
the Lie algebra
$\mathfrak{t}^{n+1}$ of $T^{n+1}$, then, as $A$ is a $T^{n+1}$-connection,
we have
$$A(\tilde Y_i)=a^\alpha_iA(X_\alpha)=a^\alpha_ie_\alpha.$$
On the other hand, we want the functions  $a^\alpha_i$ to be solutions of
the equation
$${\tilde Y_i}\rfloor dA+d({\tilde Y_i}\rfloor A)=0,$$
and hence
$$d(a^\alpha_i)=-{\tilde Y_i}\rfloor A.$$
As $Y_i$ are Hamiltonian, $Y_i\rfloor\omega^\alpha=dF^\alpha_i$, thus
$$dF^\alpha_i=Y_i\rfloor\omega^\alpha={\tilde Y_i}\rfloor A.$$
All in all $da^\alpha_i=-dF^\alpha_i$ and we may take
$a^\alpha_i=-F^\alpha_i$.

For (ii), as $T^m$ preserves $A$ and $\eta_F$ (the action of $T^m$ on $F$ is
trivial), $T^m$
preserves the contact form $\eta$ on $P\times_{T^{n+1}}F$.
\end{proof}

\begin{remark}
Suppose now that $F$ is regular and let $B_F$ the basis of its Boothby-Wang
fibration. Also,
suppose that for a torus $T$ (not necessarily of maximal dimension)
$P\times_TF$ has a
Sasakian structure. As we have seen above, this is the case when $P$ and $B$
are toric.
Lerman only constructs the contact structure on $P\times_TF$, but it can be
seen that a contact form adapted to this is written on $P\times F$ as
$\eta=f_A\cdot A +\eta_F$ for some function $f_A$ and hence
\begin{equation*}
\begin{split}
d\eta&=df_A\wedge A+f_AdA+d\eta_F\\
&=df_A\wedge A+f_A\omega_P+\omega_{B_F}
\end{split}
\end{equation*}
 Here $\omega_P$ and $\omega_{B_F}$ are $(1,1)$ forms. Splitting $df_A$ and
$A$ into
their $(0,1)$ and $(1,0)$ components, we see that $d\eta$ has a $(2,0)$
component, namely
$df_A^{(1,0)}\wedge A^{(1,0)}$, if and only if   $df_A^{(1,0)}$ and
$A^{(1,0)}$ are linearly
independent, which clearly happens when on $B\times B_F$ we take the product
complex structure. But in this case  $d\eta=\pi^*\omega$ for the K\"ahler form of
$B\times B_F$ and it has to be of type $(1,1)$. So, if the complex structure on $B\times B_F$ is
the product one, then necessarily $f_A=const.$ For $T=S^1$, this corresponds to the above
described join. Hence, the join is a particular case of Lerman's construction.
\end{remark}

\section{Toric Sasakian 5-Manifolds}

We begin this section by recalling fundamental results of Smale
and Barden concerning classification of compact smooth
simply-connected $5$-manifolds \cite{Sm62, Bar65}. Remarkably, any
such manifold is completely determined by $H_2(M,\bbz)$ and the
second Stiefel-Whitney class map $w_2$. In particular, the smooth
structure on a closed simply-connected $5$-manifold is unique.

\begin{definitiontheorem}\label{BardenInvariant}
Let $M$ be a compact, smooth, oriented, $1$-connected
$5$-manifold. Write $H_2(M,\bbz)$ as a direct sum of cyclic groups
of prime power order
\begin{equation}\label{BaSm-dec}
H_2(M,\bbz)=\bbz^k\bigoplus_{p,i} \bigl(\bbz_{p^i}\bigr)^{c(p^i)}
\end{equation}
where $k=b_2(M)$, and $c(p^i)=c(p^i,M)$. The non-negative integers
$k, c(p^i)$ are determined by $H_2(M,\bbz)$ but the subgroups
$\bbz_{p^i}\subset H_2(M,\bbz)$ are not unique. One can choose the
decomposition (\ref{BaSm-dec}) such that the second
Stiefel--Whitney class map
$$
w_2:H_2(M,\bbz)\to \bbz_2
$$ is zero
on all but one summand  $\bbz_{2^j}$. The value $j$ is unique,
denoted by $i(M)$, and called the {\bf Barden invariant} of $M$.
It can take on any value $j$ for which $c(2^j)\neq 0$, besides \
$0$ and $\infty$. Alternatively, $i(M)$ is the smallest $j$ such
that there is an $\alpha\in H_2(M,\bbz)$ such that
$w_2(\alpha)\neq 0$ and $\alpha$ has order $2^j$.
\end{definitiontheorem}

The following theorem was proved by Smale \cite{Sm62} in the spin
case in when $w_2=0$ implying  $i=0$. Subsequent generalization
with no assumption on $w_2$ is due to Barden \cite{Bar65}. We
shall formulate it here using Barden's notation.

\begin{theorem}\label{Barden-Smale-Theorem}
The class $\calb$ of simply connected, closed, oriented, smooth,
5-manifolds is classifiable under diffeomorphism. Furthermore, any
such $M$ is diffeomorphic to one of the spaces
\begin{equation}\label{Barden-decomposition}
M_{j;k_1,\ldots,k_s}=X_j\#M_{k_1}\#\cdots \#M_{k_s},
\end{equation} where $-1\leq j\leq\infty,$ $s\geq0$, $1<k_1$ and
$k_i$ divides $k_{i+1}$ or $k_{i+1}=\infty$. A complete set of
invariants is provided by $H_2(M,\bbz)$ and $i(M)$ and the
manifolds $X_{-1},X_0,X_j,X_\infty,M_j,M_\infty$ are characterized
as follows
\begin{center}
\begin{tabular}{|l|c|c|}\hline
$M$ & $H_2(M,\bbz)$ &$i(M)$\\
\hline\hline
$X_{-1}=SU(3)/SO(3)$&$\bbz_2$&$1$ \\
\hline
$M_0=X_0=S^5$ & $0$ &$0$ \\
\hline
$X_j$, \ \ $0<j<\infty$ &$\bbz_{2^j}\oplus\bbz_{2^j}$ &$j$ \\
\hline
$X_\infty$&$\bbz$&$\infty$ \\
\hline
$M_k$, \ \ $0<k<\infty$&$\bbz_k\oplus\bbz_k$&$0$ \\
\hline
$M_\infty=S^2\times S^3$&$\bbz$&$0$ \\
\hline
\end{tabular}
\end{center}
\end{theorem}
\bigskip

In this section we would like to investigate the question which of
the manifolds in $\calb$ admit toric Sasakian structures. We begin
with an important example.

\begin{example} {\bf [Circle Bundles over Hirzebruch surfaces]}
We shall construct infinite families of deformation classes of
toric Sasakian structures on circle bundles over Hirzebruch
surfaces. The total space $M$ of any such circle bundle must have
$b_1(M)=0$ and $b_2(M)=1.$ Theorem \ref{Barden-Smale-Theorem}
implies that there are precisely two 5-manifolds  in $\calb$ with
$b_2=1.$ They are $M_\infty=S^2\times S^3$ and $X_\infty$, the
non-trivial $S^3$-bundle over $S^2$. They both have
$H_2(M,\bbz)=\bbz$ and are distinguished by their Barden
invariant.

Recall the Hirzebruch surfaces $S_n$ (cf. \cite{GrHa78}, pgs 517-520) are realized as
the projectivizations of the sum of two line bundles over
$\bbc\bbp^1,$ which we can take as
$$S_n=\bbp\bigl(\calo \oplus\calo(n)\bigr).$$
They are diffeomorphic to $\bbc\bbp^1\times \bbc\bbp^1$ if $n$ is even, and to the
blow-up of $\bbc\bbp^2$ at one point, which we denote as $\widetilde{\bbc\bbp}^2,$ if
$n$
is odd. For $n=0$ and $1$ we get $\bbc\bbp^1\times \bbc\bbp^1$ and
$\widetilde{\bbc\bbp}^2,$ respectively. Now ${\rm Pic}(S_n)\approx \bbz\oplus \bbz,$
and we can take the
Poincar\'e duals of a section of $\calo(n)$ and the homology class of the fibre as its
generators. The corresponding divisors can be represented by rational curves which we
denote by $C$ and $F,$ respectively satisfying
$$C\cdot C=n, \qquad F\cdot F=0,\qquad C\cdot F=1.$$
Let $\gra_1$ and $\gra_2$ denote the Poincar\'e duals of $C$ and $F$, respectively. The
classes
$\gra_1$ and $\gra_2$ can be represented by $(1,1)$ forms $\gro_1$ and $\gro_2,$
respectively,
so that the $(1,1)$ form $\gro_{l_1,l_2}=l_1\gro_1 +l_2\gro_2$ determines a circle bundle
over
$S_n$ whose first Chern class is $[\gro_{l_1,l_2}].$ We thus have circle bundles
depending on a triple of integers $(l_1,l_2,n),$ with $n$ non-negative,
$$S^1\ra{1.7}M_{l_1,l_2,n}\fract{\pi}{\ra{1.7}}S_n.$$
Now in order that $M_{l_1,l_2,n}$ admit a Sasakian structure it is necessary that
$\gro_{l_1,l_2}$ be a
positive $(1,1)$ form, that is, $\gro_{l_1,l_2}$ must lie within the K\"ahler cone
$\calk(S_n).$ The conditions for positivity are by Nakai's criterion,
\begin{enumerate}
\item $\gro_{l_1,l_2}^2>0,$
\item $\int_D\gro_{l_1,l_2} > 0$ for all holomorphic curves $D,$
\end{enumerate}
which in our case give $l_1,l_2>0.$

Next we determine the diffeomorphism type of $M_{l_1,l_2,n}.$ Since the K\"ahler class
$[\gro_{l_1,l_2}]$ transgresses to the derivative of the contact form, $d\eta_{l_1,l_2},$ we
see
that $\pi^*\gra_1 =-l_2\grg$ and $\pi^*\gra_2=l_1\grg$ where $\grg$ is a generator of
$H^2(M_{l_1,l_2,n},\bbz)\approx \bbz.$ Now the first Chern class of $S_n$ is \cite{GrHa78}
\begin{equation}\label{HirChern}
c_1(S_n)=2\gra_1-(n-2)\gra_2,
\end{equation}
which pulls back to the basic first Chern class on
$M_{l_1,l_2,n}.$ So the first Chern class of the contact line
bundle $\cald$ is given by
\begin{equation}\label{HirChern2}
c_1(\cald)=-[2l_2-l_1(2-n)]\grg.
\end{equation}
If we take the integers $l_1,l_2$ to be relatively prime then the
manifold $M_{l_1,l_2,n}$ will be simply connected. Furthermore,
since $M_{l_1,l_2,n}$ has a regular contact structure, there is no
torsion in $H_2(M_{l_1,l_2,n},\bbz)$ \cite{Gei91}. Thus, by
Theorem \ref{Barden-Smale-Theorem} $M_{l_1,l_2,n}$ is either
$S^2\times S^3$ or $X_\infty,$ depending on whether $M_{l_1,l_2,n}$ is
spin or not. But we have $w_2(M_{l_1,l_2,n})\equiv nl_1 \mod~2,$
so $M_{l_1,l_2,n}$ is diffeomorphic to $S^2\times S^3$ if $nl_1$
is even, and to $X_\infty$ if $nl_1$ is odd. Now that we have
identified the manifolds $M_{l_1,l_2,n}$, we wish to distinguish
the deformation classes of Sasakian structures that live on them.
Thus, we consider the equivalence classes of regular homologous
Sasakian structures $\gF_{l_1,l_2,n}$ \cite{BGN03b,BG05}.
We have arrived at

\begin{theorem}\label{circlebunHir1}
For each triple of positive integers $(l_1,l_2,m)$ satisfying $\gcd(l_1,l_2)=1,$ the manifold
$S^2\times S^3$ admits the following deformation classes of regular Sasakian structures
$\gF_{l_1,l_2,2m}$ and $\gF_{2l_1,l_2,2m+1}.$
\end{theorem}

\begin{theorem}\label{circlebunHir2}
For each triple of positive integers $(l_1,l_2,m)$ satisfying $\gcd(l_1,l_2)=1,$ the manifold
$X_\infty$ admits the deformation classes of regular Sasakian structures
$\gF_{2l_1-1,l_2,2m-1}.$
\end{theorem}

Now we have

\begin{theorem}\label{toricHir}
The deformation classes of Sasakian structures described in Theorems \ref{circlebunHir1}
and
\ref{circlebunHir2} are all toric.
\end{theorem}

\begin{proof}
As in \cite{Aud94} we describe Hirzebruch surfaces $S_n$ as smooth algebraic
subvarieties
of $\bbc\bbp^1\times \bbc\bbp^2$, viz.
$$S_n =\{[u_0,u_1],[v_0,v_1,v_2] ~|~u_0^nv_1=u_1^nv_0\}\subset \bbc\bbp^1\times
\bbc\bbp^2.$$
Now $S_n$ admits the action of a complex 2-torus by
$$([u_0,u_1],[v_0,v_1,v_2])\mapsto ([\grt u_0,\grz
u_1],[\grz^{-n}v_0,\grt^{-n}v_1,\grt^{-n}\grz^{-n}v_2]).$$
The homology class $F$ is represented by the rational curve $([a,b],[0,0,1]),$ while $C$ is
represented by $([a,b],[a^n,b^n,0]),$ and it is easy to check that these rational curves are
invariant under $\bbc^*\times \bbc^*$-action given above. It follows that for each
admissible
value of $(l_1,l_2)$ the K\"ahler form $\gro=l_1\gro_1+l_2\gro_2$ is invariant under the
toral
subgroup $\gT_2$ of $\bbc^*\times \bbc^*$. Furthermore, the action of $\gT_2$ is
Hamiltonian and hence, it lifts to a $\gT_2$ in the automorphism group of the Sasakian
structure. This together with the $S^1$ generated by the Reeb vector field $\xi_{l_1,l_2}$
gives $M_{l_1,l_2,n}$ a toric Sasakian  structure.
\end{proof}

\begin{remark}\label{Hirzremark}
The existence of toric Sasakian structures on $S^2\times S^3$ and
$X_\infty$ also follows from Theorem \ref{Sasfibre} as well as
Theorem \ref{toricSas5man} below.
\end{remark}

Next we briefly discuss which toric Sasakian structures $\gF_{l_1,l_2,n}$ belong to equivalent 
contact structures.
It is convenient to make a change of basis of $H^2(S_n,\bbz).$ For simplicity we consider the case 
$n=2m$ so
$S_{2m}$ is diffeomorphic to $S^2\times S^2.$ For $i=1,2$ we let $\grs_i$ denote the classes in 
$H^2(S^2\times
S^2,\bbz)$ given by pulling back the volume form on the $i^{th}$ factor. Writing the K\"ahler class 
$[\gro]=
a_1\grs_1+a_2\grs_2$ in terms of the this basis, we see that
$$a_1=l_1m+l_2, \qquad a_2=l_1,$$
and the positivity condition becomes $a_1> ma_2>0.$ We denote the corresponding deformation 
classes of toric
Sasakian structures on $S^2\times S^3$ by $\gF(a_1,a_2,m).$ The integers $a_1,a_2$ are written 
as $a,b$ in
\cite{Kar03} and \cite{Ler03b}. In terms of the $a_i$ the first Chern class \ref{HirChern2} simplifies 
to $c_1(\cald)=
2(a_1-a_2)\grg.$ Thus, $\gF(a_1,a_2,m)$ and $\gF(a'_1,a'_2,m')$ belong to non-isomorphic 
contact structures if
$a'_1-a'_2\neq a_1-a_2.$ The following theorem is due to Lerman \cite{Ler03b}.

\begin{theorem}\label{circlebunHir3}
For every pair of relatively prime integers $(a_1,a_2)$ there are
$\lceil \frac{a_2}{a_1}\rceil$ inequivalent regular toric Sasakian
structures on $S^2\times S^3$ having the same contact form
$\eta_{a_1,a_2}.$ However, for each integer $m=0,\cdots, \lceil \frac{a_2}{a_1}\rceil$ the structures 
$\gF(a_1,a_2,m)$
are inequivalent as toric contact structures.
\end{theorem}

Here $\lceil a \rceil$ denotes the smallest integer greater than
or equal to $a.$ The essence of Theorem \ref{circlebunHir3} is that $\lceil\frac{a_2}{a_1}\rceil$ is
precisely the number of non-conjugate maximal tori in the
contactomorphism group $\gC\go\gn(S^2\times S^3,\eta_{a_1,a_2})$ \cite{Kar03,Ler03b}.
\end{example}

We now begin the discussion of the general toric case. It turns
out that even asking for just an effective $T^3$ action on $M$
severely restricts its topology. Recall the following
classification theorem of Oh \cite{Oh83}:

\begin{theorem}\label{Ohthm}
Let $M$ be a closed simply connected 5-manifold with an effective
$T^3$-action. Then $M$ has no 2-torsion. In particular, $M$ is
diffeomorphic to $S^5, k(S^2\times S^3)$ or $X_\infty\#(k-1)(S^2\times S^3)$,
where $k=b_2(M)\geq 1$. Conversely, all these manifolds admit
effective $T^3$ actions.
\end{theorem}

In particular, there are infinitely many Sasakian 5-manifolds
(even Sasakian-Einstein) which do not admit any toric contact
structure. In \cite{Yam01} Yamazaki proved that all Oh's toric
5-manifolds also admit compatible K-contact structures. But then
the main theorem in \cite{BG00b} can be used to strengthen this to

\begin{theorem}\label{toricSas5man}
Let $M$ be a closed simply connected 5-manifold with an effective
$T^3$-action. Then $M$ admits toric Sasakian structures and is
diffeomorphic to $S^5, k(S^2\times S^3)$ or $X_\infty\#(k-1)(S^2\times S^3)$,
where $k=b_2(M)\geq 1.$
\end{theorem}

Further recall that Geiges \cite{Gei91} showed that the torsion in
$H_2(M^5,\bbz)$ is the only obstruction to the existence of a {\it
regular} contact structure on $M$. So the question arises whether,
for a given $M$ in Theorem \ref{toricSas5man} there exist a
regular Sasakian structures compatible with some toric contact
structure. We answer this in the affirmative by giving an explicit
construction as circle bundles over the blow-ups of the Hirzebruch
surfaces.

It is well-known \cite{Ful93} that the smooth toric surfaces are
all obtained by blowing-up Hirzebruch surfaces at the fixed points
of the $T^2$ action. Begin with a Hirzebruch surface $S_n$ and
blow-up $S_n$ at one of the 4 fixed points of the $T^2$ action.
This gives a smooth toric surface $S_{n,1}$ which can be
represented by a Delzant polytope with 5 vertices. Repeat this
procedure inductively to obtain smooth algebraic toric surfaces
$S_{n,k}$ whose Delzant polytope has $k+4$ vertices. Choose a
K\"ahler class $[\gro]$ lying on the Neron-Severi lattice, and
construct the circle bundle $\pi_{n,k}:M_{n,k}\ra{1.3} S_{n,k}$
whose Euler class is $[\gro].$ The K\"ahler form can be chosen to
be invariant under the $T^2$ action, and we can choose a $T^2$
invariant connection $\eta$ in $\pi_{n,k}:M_{n,k}\ra{1.3} S_{n,k}$
whose curvature form satisfies $d\eta=\pi_{n,k}^*\gro.$ This gives
a regular Sasakian structure $\cals=(\xi,\eta,\Phi,g)$ on
$M_{n,k},$ and by Theorem \ref{Barden-Smale-Theorem} $M_{n,k}$ is
diffeomorphic to either $(k+1)(S^2\times S^3)$ or $X_\infty\#k (S^2\times S^3).$
Since $H^1(S_{n,k},\bbz)=0$ the torus action $T^2$ lifts to a
$T^2$ action in the automorphism group of the Sasakian structure
$\cals$ (cf. \cite{BG05}) which together with the circle group
generated by the Reeb field $\xi$ makes $\cals$ a regular toric
Sasakian structure.

It remains to show that all of the manifolds in $\calb$ with no
2-torsion occur. For this we need to compute the second
Stiefel-Whitney class $w_2(M_{n,k})$ which is the mod 2 reduction
of the first Chern class of contact bundle $\cald_{n,k}$ of
$M_{n,k}.$ First we give a K\"ahler form of the complex manifolds
$S_{n,k}$ constructed in the paragraph above. Let
$\tilde{\gro}_{l_1,l_2}, \tilde{\gro}_1,\tilde{\gro}_2$ denote the
proper transform of $\gro_{l_1,l_2},\gro_1,\gro_2$ respectively.
Then the K\"ahler form on $S_{n,k}$ can be written as
\begin{equation}\label{Kahblowup}
\tilde{\gro}_{l_1,\cdots,l_{k+2}} = \sum_il_i \tilde{\gro}_i
\end{equation}
where $\tilde{\gro}_{i+2}$ is the $(1,1)$ form representing the
Poincar\'e duals $\tilde{\gra_{i+2}}$ of the exceptional divisors
$E_i.$  Then writing the K\"ahler class as
$\sum_{i=1}^kl_{i}\tilde{\gra}_{i}$, the positivity condition
becomes
\begin{equation}\label{posblowup}
0<\sum_{i=1}^kl_{i}\tilde{\gra}_{i}\cup
\sum_{i=1}^kl_{i}\tilde{\gra}_{i}=l_1(2l_2+nl_1)-\sum_{i=1}^kl_{i+1}^2.
\end{equation}
So $\tilde{\gro}_{l_1,\cdots,l_{k+2}}$ defines a K\"ahler metric
if this inequality is satisfied. Here we have used the fact the
exceptional divisor is in the kernel of the corresponding blow-up
map. Define the integer valued $k+2$-vector
$\bfl=(l_1,\cdots,l_{k+2}).$ It is convenient to choose
$\bfl=(1,l_2,1,\cdots,1)$ in which case the positivity condition
becomes $2l_2+n>k.$ For simplicity we denote the corresponding
K\"ahler form by $\tilde{\gro}_{l_2}.$

Next we need to compute the first Chern class of $\cald_{n,k}.$ This is
$\pi_{n,k}^*c_1(S_{n,k})$ modulo the transgression of  the K\"ahler class on $S_{n,k},$ that
is, modulo the relation $\pi_{n,k}^*\sum_{i=1}^{k+2}l_i\tilde{\gra}_i =0.$  Let
$\grb_1,\cdots,\grb_{k+1}$ be a basis for $H^2(M_{n,k},\bbz),$ and write
$\pi_{n,k}^*\tilde{\gra_i}=\sum_{j=1}^{k+1}m_{ij}\grb_j.$ We need to choose the $k+2$ by
$k+1$ matrix $(m_{ij})$ such that \begin{equation}\label{mlortho}
\sum_{i=1}^{k+2}l_im_{ij}=0.
\end{equation}
Now using Equation \ref{HirChern} we have
\begin{equation}\label{c1blowup}
c_1(S_{n,k})=2\tilde{\gra}_1-(n-2)\tilde{\gra}_2-\sum_{i=1}^k\tilde{\gra}_{k+2}
\end{equation}
which gives
\begin{equation}\label{c1blowup2}
\pi_{n,k}^*c_1(S_{n,k})=2\sum_jm_{1j}\grb_j-(n-2)\sum_jm_{2j}\grb_j-\sum_{i=3}^{k+2}
\sum_jm_{ij}\grb_j.
\end{equation}
We now make a judicious choice of the matrix $(m_{ij}).$
\begin{equation}\label{mmatrix}
(m_{ij}) =\left(
\begin{matrix}
-l_2 & 2 & 2 & \cdots & 2 \\
1  & 0 & 0 & \cdots & 0 \\
0 & -2 & \ddots & \vdots & \vdots \\
\vdots & 0 & \ddots & 0 & 0 \\
\vdots& \vdots & \cdots & -2 & 0 \\
0     & 0 & \cdots & 0 & -2
\end{matrix}
\right)
\end{equation}
The orthogonality condition \ref{mlortho} is satisfied and Equation \ref{c1blowup} becomes
$$\pi_{n,k}^*c_1(S_{n,k})= 2l_2\grb_1 -6(\grb_2+\cdots \grb_{k+1}) -(n-2)\grb_1.$$
It follows that $w_2(M_{n,k})\equiv n\mod 2.$

We have arrived at
\begin{theorem}\label{toricSas5man2}
Let $S_{n,k}$ be the equivariant $k$-fold blow-up of the
Hirzebruch surface $S_n.$ Let $\pi_{n,k}:M_{n,k}\ra{1.5} S_{n,k}$
be the circle bundle defined by the integral K\"ahler form
$\tilde{\gro}_{l_2}.$ Then for each positive integer $l_2$ satisfying
$2l_2+n>k,$ the manifold $M_{n,k}$ admits a toric regular Sasakian
structure, $M_{n,k}$ is diffeomorphic to $k(S^2\times S^3)$ if $n$ is
even, and if $n$ is odd it is diffeomorphic to
$X_\infty\#(k-1)(S^2\times S^3).$ Thus, every regular contact 5-manifold
admits a toric regular Sasakian structure.
\end{theorem}

\begin{remark}
With more analysis one can make a count of inequivalent
deformation classes of toric Sasakian structures on the manifolds
$k(S^2\times S^3)$ and $X_\infty\#(k-1)(S^2\times S^3)$ as done in Theorems
\ref{circlebunHir1} and \ref{circlebunHir2}.

We close this section with a brief discussion of toric
Sasakian-Einstein structures. It is well-known that $S^5$ and
$S^2\times S^3$ have homogeneous, hence regular, Sasakian-Einstein
structures. Both of these are clearly toric. First inhomogeneous
examples of explicit toric Sasakian-Einstein structures on
$S^2\times S^3$ were obtained by Gauntlett et al. \cite{GMSW04a,
MaSp06}. These metrics are of cohomogeneity 1 with $U(2)\times
U(1)$ acting by isometries. In fact, Gauntlett et al. construct
infinite families of toric Sasakian-Einstein structures
parameterized by two relatively prime integers $p>q$. When
$4p^2-3q^2=n^2$ their examples are quasi-regular, i.e., the Reeb
vector field has closed orbits. Otherwise the Sasakian-Einstein
structure is not quasi-regular. These new examples were further
generalized by Cveti\v c et al. \cite{CLPP05} (see also
\cite{MaSp05b}) who found toric Sasakian-Einstein metrics on
$S^2\times S^3$ of cohomogeneity 2.

This raises a natural question: Do all spin manifolds of Theorem
\ref{toricSas5man} admit a toric Sasakian-Einstein structure?
Since simply connected Sasakian-Einstein spaces are necessarily
spin, $X_\infty\#(k-1)(S^2\times S^3)$ must be excluded. It is known that
$S^5$ and $k(S^2\times S^3)$ admit families of quasi-regular
Sasakian-Einstein structures for any $k$ \cite{BGN03c, BGK05,
Kol04}. But most of these metrics have only a one-dimensional
isometry group. However, it turns out that $S^2\times S^3$ is by
no means special in this respect: there exist families of toric
Sasakian-Einstein structures on $k(S^2\times S^3)$ for arbitrary $k$. This has just recently been 
proven 
in \cite{CFO07,FOW06}, and a bit earlier van Coevering \cite{Coe06} proved this result for $k$ odd.
\end{remark}

\def\cprime{$'$} \def\cprime{$'$} \def\cprime{$'$} \def\cprime{$'$}
\providecommand{\bysame}{\leavevmode\hbox to3em{\hrulefill}\thinspace}
\providecommand{\MR}{\relax\ifhmode\unskip\space\fi MR }
\providecommand{\MRhref}[2]{%
  \href{http://www.ams.org/mathscinet-getitem?mr=#1}{#2}
}
\providecommand{\href}[2]{#2}


\begin{thebibliography}{GMSW04}

\bibitem[Aud94]{Aud94}
M.~Audin, \emph{Symplectic and almost complex manifolds}, Holomorphic curves in
  symplectic geometry, Progr. Math., vol. 117, Birkh\"auser, Basel, 1994, With
  an appendix by P. Gauduchon, pp.~41--74. \MR{1274926}

\bibitem[Bar65]{Bar65}
D.~Barden, \emph{Simply connected five-manifolds}, Ann. of Math. (2)
  \textbf{82} (1965), 365--385. \MR{32 \#1714}

\bibitem[BG00a]{BG00b}
C.~P. Boyer and K.~Galicki, \emph{A note on toric contact geometry}, J. Geom.
  Phys. \textbf{35} (2000), no.~4, 288--298. \MR{2001h:53124}

\bibitem[BG00b]{BG00a}
\bysame, \emph{On {S}asakian-{E}instein geometry}, Internat. J. Math.
  \textbf{11} (2000), no.~7, 873--909. \MR{2001k:53081}

\bibitem[BG03]{BG03}
Charles~P. Boyer and Krzysztof Galicki, \emph{New {E}instein metrics on
  {$8\#(S^2\times S^3)$}}, Differential Geom. Appl. \textbf{19} (2003), no.~2,
  245--251. \MR{2 002 662}

\bibitem[BG06]{BG05}
C.~P. Boyer and K.~Galicki, \emph{{Sasakian Geometry}}, Oxford Mathematical
  Monographs, Oxford University Press, to appear, Oxford, 2006.

\bibitem[BGK05]{BGK05}
C.~P. Boyer, K.~Galicki, and J.~Koll\'ar, \emph{Einstein metrics on spheres},
  Ann. of Math. (2) \textbf{162} (2005), no.~1, 557--580.

\bibitem[BGKT05]{BGKT05}
C.~P. Boyer, K.~Galicki, J.~Koll{\'a}r, and E.~Thomas, \emph{Einstein metrics
  on exotic spheres in dimensions 7, 11, and 15}, Experiment. Math. \textbf{14}
  (2005), no.~1, 59--64. \MR{2146519}

\bibitem[BGM06]{BGM06}
C.~P. Boyer, K.~Galicki, and P.~Matzeu, \emph{{On Eta-Einstein Sasakian
  Geometry}}, Comm. Math. Phys. \textbf{262} (2006), 177--208.

\bibitem[BGN02]{BGN02b}
C.~P. Boyer, K.~Galicki, and M.~Nakamaye, \emph{Sasakian-{E}instein structures
  on {$9\#(S^2\times S^3)$}}, Trans. Amer. Math. Soc. \textbf{354} (2002),
  no.~8, 2983--2996 (electronic). \MR{2003g:53061}

\bibitem[BGN03a]{BGN03c}
\bysame, \emph{On the geometry of {S}asakian-{E}instein 5-manifolds}, Math.
  Ann. \textbf{325} (2003), no.~3, 485--524. \MR{2004b:53061}

\bibitem[BGN03b]{BGN03b}
\bysame, \emph{Sasakian geometry, homotopy spheres and positive {R}icci
  curvature}, Topology \textbf{42} (2003), no.~5, 981--1002. \MR{1 978 045}

\bibitem[BGN03c]{BGN03a}
Charles~P. Boyer, Krzysztof Galicki, and Michael Nakamaye, \emph{On positive
  {S}asakian geometry}, Geom. Dedicata \textbf{101} (2003), 93--102.
  \MR{2017897 (2005a:53072)}

\bibitem[BM93]{BM93}
A.~Banyaga and P.~Molino, \emph{G\'eom\'etrie des formes de contact
  compl\`etement int\'egrables de type toriques}, S\'eminaire Gaston Darboux de
  G\'eom\'etrie et Topologie Diff\'erentielle, 1991--1992 (Montpellier), Univ.
  Montpellier II, Montpellier, 1993, pp.~1--25. \MR{94e:53029}

\bibitem[CFO07]{CFO07}
K. Cho, A. Futaki, and H. Ono, \emph{{Uniqueness and examples of
  compact toric {S}asaki-{E}instein metrics}}, preprint, arXiv:math.DG/0701122 (2007).

\bibitem[CLPP05]{CLPP05}
M.~Cveti{\v{c}}, H.~L{\"u}, Don~N. Page, and C.~N. Pope, \emph{New
  {E}instein-{S}asaki spaces in five and higher dimensions}, Phys. Rev. Lett.
  \textbf{95} (2005), no.~7, 071101, 4. \MR{2167018}

\bibitem[DO98]{DrOr98}
S.~Dragomir and L.~Ornea, \emph{Locally conformal {K}\"ahler geometry},
  Progress in Mathematics, vol. 155, Birkh\"auser Boston Inc., Boston, MA,
  1998. \MR{1481969 (99a:53081)}

\bibitem[Ful93]{Ful93}
W.~Fulton, \emph{Introduction to toric varieties}, Annals of Mathematics
  Studies, vol. 131, Princeton University Press, Princeton, NJ, 1993, The
  William H. Roever Lectures in Geometry. \MR{94g:14028}

\bibitem[FOW06]{FOW06}
A.~Futaki, H.~Ono, and G.~Wang, \emph{{Transverse K\"ahler geometry of Sasaki
  manifolds and toric Sasaki-Einstein manifolds}}, preprint: arXiv:math.DG/0607586 (2006).

\bibitem[Gei91]{Gei91}
H.~Geiges, \emph{Contact structures on {$1$}-connected {$5$}-manifolds},
  Mathematika \textbf{38} (1991), no.~2, 303--311 (1992). \MR{93e:57042}

\bibitem[GH78]{GrHa78}
P.~Griffiths and J.~Harris, \emph{Principles of algebraic geometry},
  Wiley-Interscience [John Wiley \& Sons], New York, 1978, Pure and Applied
  Mathematics. \MR{80b:14001}

\bibitem[GK05]{GhKo05}
A~Ghigi and J.~Koll\'ar, \emph{{K\"ahler Einstein metrics on orbifolds and
  Einstein metrics on Spheres}}, arXiv:math.DG/0507289 (2005).

\bibitem[GLS96]{GuLeSt96}
V.~Guillemin, E.~Lerman, and S.~Sternberg, \emph{Symplectic fibrations and
  multiplicity diagrams}, Cambridge University Press, Cambridge, 1996.
  \MR{1414677 (98d:58074)}

\bibitem[GMSW04]{GMSW04a}
J.~P. Gauntlett, D.~Martelli, J.~Sparks, and D.~Waldram,
  \emph{Sasaki-{E}instein metrics on {$S^2\times S^3$}}, Adv. Theor. Math.
  Phys. \textbf{8} (2004), no.~4, 711--734. \MR{2141499}

\bibitem[Hae84]{Hae84}
A.~Haefliger, \emph{Groupo\"\i des d'holonomie et classifiants}, Ast\'erisque
  (1984), no.~116, 70--97, Transversal structure of foliations (Toulouse,
  1982). \MR{86c:57026a}

\bibitem[HS91]{HaSa91}
A.~Haefliger and {\'E}.~Salem, \emph{Actions of tori on orbifolds}, Ann. Global
  Anal. Geom. \textbf{9} (1991), no.~1, 37--59. \MR{92f:57047}

\bibitem[Kar03]{Kar03}
Y.~Karshon, \emph{Maximal tori in the symplectomorphism groups of {H}irzebruch
  surfaces}, Math. Res. Lett. \textbf{10} (2003), no.~1, 125--132.
  \MR{1960129 (2004f:53101)}

\bibitem[Kol04]{Kol04}
J.~Koll\'ar, \emph{{Einstein metrics on connected sums of $S^2\times S^3$}},
  arXiv:math.DG/0402141 (2004).

\bibitem[KS88]{KS88}
M.~Kreck and S.~Stolz, \emph{A diffeomorphism classification of
  {$7$}-dimensional homogeneous {E}instein manifolds with {${\rm
  SU}(3)\times{\rm SU}(2)\times{\rm U}(1)$}-symmetry}, Ann. of Math. (2)
  \textbf{127} (1988), no.~2, 373--388. \MR{89c:57042}

\bibitem[Ler02]{Ler02a}
E.~Lerman, \emph{Contact toric manifolds}, J. Symplectic Geom. \textbf{1}
  (2002), no.~4, 785--828. \MR{2 039 164}

\bibitem[Ler03]{Ler03b}
\bysame, \emph{Maximal tori in the contactomorphism groups of circle bundles
  over {H}irzebruch surfaces}, Math. Res. Lett. \textbf{10} (2003), no.~1,
  133--144. \MR{2004g:53097}

\bibitem[Ler04a]{Ler04b}
\bysame, \emph{Contact fiber bundles}, J. Geom. Phys. \textbf{49} (2004),
  no.~1, 52--66. \MR{2077244}

\bibitem[Ler04b]{Ler04}
\bysame, \emph{Homotopy groups of {$K$}-contact toric manifolds}, Trans. Amer.
  Math. Soc. \textbf{356} (2004), no.~10, 4075--4083 (electronic). \MR{2 058
  839}

\bibitem[Mol88]{Mol88}
P.~Molino, \emph{Riemannian foliations}, Progress in Mathematics, vol.~73,
  Birkh\"auser Boston Inc., Boston, MA, 1988, Translated from the French by
  Grant Cairns, With appendices by Cairns, Y. Carri\`ere, \'E. Ghys, E. Salem
  and V. Sergiescu. \MR{89b:53054}

\bibitem[MS05]{MaSp05b}
D.~Martelli and J.~Sparks, \emph{Toric {S}asaki-{E}instein metrics on {$S\sp
  2\times S\sp 3$}}, Phys. Lett. B \textbf{621} (2005), no.~1-2, 208--212.
  \MR{2152673}

\bibitem[MS06]{MaSp06}
\bysame, \emph{Toric geometry, {S}asaki-{E}instein manifolds and a new infinite
  class of {A}d{S}/{C}{F}{T} duals}, Comm. Math. Phys. \textbf{262} (2006),
  51--89.

\bibitem[Oh83]{Oh83}
H.~S. Oh, \emph{Toral actions on {$5$}-manifolds}, Trans. Amer. Math. Soc.
  \textbf{278} (1983), no.~1, 233--252. \MR{697072 (85b:57043)}

\bibitem[OV03]{OrVe03}
L.~Ornea and M.~Verbitsky, \emph{Structure theorem for compact {V}aisman
  manifolds}, Math. Res. Lett. \textbf{10} (2003), no.~5-6, 799--805.
  \MR{2024735 (2004j:53093)}

\bibitem[OV05]{OrVe05}
\bysame, \emph{An immersion theorem for {V}aisman manifolds}, Math. Ann.
  \textbf{332} (2005), no.~1, 121--143. \MR{2139254}

\bibitem[Sma62]{Sm62}
S.~Smale, \emph{On the structure of {$5$}-manifolds}, Ann. of Math. (2)
  \textbf{75} (1962), 38--46. \MR{25 \#4544}

\bibitem[Ste51]{Stee51}
N.~Steenrod, \emph{The {T}opology of {F}ibre {B}undles}, Princeton Mathematical
  Series, vol. 14, Princeton University Press, Princeton, N. J., 1951.
  \MR{12,522b}

\bibitem[vC]{Coe06}
C.~van Coevering, \emph{{Toric Surfaces and Sasakian-Einstein 5-manifolds}},
  SUNY at Stony Brook Ph.D. Thesis. preprint, arXiv:math.DG/0607721 (2006).

\bibitem[Wad75]{Wad}
A.~W. Wadsley, \emph{Geodesic foliations by circles}, J. Differential Geometry
  \textbf{10} (1975), no.~4, 541--549. \MR{53 \#4092}

\bibitem[Wei80]{Wei80}
A.~Weinstein, \emph{Fat bundles and symplectic manifolds}, Adv. in Math.
  \textbf{37} (1980), no.~3, 239--250. \MR{82a:53038}

\bibitem[WZ90]{WaZi90}
M.~Y. Wang and W.~Ziller, \emph{Einstein metrics on principal torus bundles},
  J. Differential Geom. \textbf{31} (1990), no.~1, 215--248. \MR{91f:53041}

\bibitem[Yam99]{Yam99}
T.~Yamazaki, \emph{A construction of {$K$}-contact manifolds by a fiber join},
  Tohoku Math. J. (2) \textbf{51} (1999), no.~4, 433--446. \MR{2001e:53094}

\bibitem[Yam01]{Yam01}
\bysame, \emph{On a surgery of {$K$}-contact manifolds}, Kodai Math. J.
  \textbf{24} (2001), no.~2, 214--225. \MR{2002c:57048}

\end{thebibliography}
\end{document}